\documentclass{amsart}
\usepackage[english]{babel}
\usepackage{amssymb,textcomp,bbm,latexsym}
\usepackage{amsrefs}
\usepackage{amsthm}
\usepackage{color}
\newcommand{\assign}{:=}
\newcommand{\mathd}{\mathrm{d}}
\newcommand{\tmop}[1]{\ensuremath{\operatorname{#1}}}

\theoremstyle{plain} 

\newtheorem{theorem}[equation]{Theorem}

\theoremstyle{definition}
\newtheorem{definition}[equation]{Definition} 

\theoremstyle{remark}

\numberwithin{equation}{section}

\begin{document}

\title{Failure of the matrix weighted bilinear Carleson embedding theorem}

\author{Komla Domelevo}

\author{Stefanie Petermichl}{\thanks{work partially supported by ERC grant CHRiSHarMa no. DLV-682402}}

\author{Kristina Ana \v{S}kreb}

%
%
%
%

%
\begin{abstract}
We prove failure of the natural formulation of a matrix weighted bilinear
    Carleson embedding theorem, featuring a matrix valued Carleson sequence as well as
    products of norms for the embedding. We show that assuming an $A_2$
    weight is also not sufficient. Indeed, a uniform bound on the conditioning
    number of the matrix weight is necessary and sufficient to get the
    bilinear embedding. We prove the optimal dependence of the embedding on
    this quantity. We show that any improvement of a recent matrix weighted
    bilinear embedding, featuring a scalar Carleson sequence and inner
    products instead of norms must fail. In particular, replacing the scalar
    sequence by a matrix sequence results in failure 
    even when maintaining the formulation using inner products. Any
    formulation using norms, even in the presence of a scalar Carleson
    sequence must fail. As a positive result, we prove the so--called matrix
    weighted redundancy condition in full generality.
\end{abstract}

{\maketitle}

\

\section{Introduction}

The Carleson embedding theorem (CET) is a classical theorem in harmonic
analysis with many applications to PDE. It states that a Carleson measure
gives an $L^2$ embedding for a function. The Carleson embedding theorem first appeared in L. Carleson's solution of the free interpolation problem \cite{Ca1958} 
and later was used in his celebrated proof of the Corona Theorem \cite{Ca1962}. 

In this paper we are primarily using the language of dyadic cubes. See for example \cite{NT1997} for the formulation therein and the illustration of its proof via Bellman functions.
One may consider any underlying measure here, so the weighted version is just the same as the unweighted
version with a very similar proof. 

A weighted bilinear embedding theorem (BET) was an important tool in
the early days of sharp weighted theory. Indeed, it was a crucial ingredient
in the first sharp weighted estimates for classical singular operators \cite{P2007}. It states that a Carleson measure
gives rise to a bilinear estimate, featuring two different functions. In this
note we are concerned with the failure of its matrix weighted analogs.

\

An unweighted CET with matrix Carleson measure holds trivially, derived from
the scalar case. Other than in the scalar case, the extension to the weighted
setting is not trivial. First versions imposed the so--called $A_2$ property
of the matrix weight introduced by Treil--Volberg in \cite{TV1997}, a condition that had been absent in the scalar case.
Recently Culiuc--Treil \cite{CT} obtained the matrix weighted version of this
theorem without any restriction on the weight, other than it being a matrix
weight.

\

The so--called matrix $A_2$ conjecture asks for the exact
growth of the matrix weighted norm estimate of the Hilbert transform acting on vector
functions. The conjectured growth estimate is that of a linear dependence on the matrix $A_2$ characteristic of the weight. 
Motivated by this problem, we consider the question of
a bilinear version of the matrix weighted Carleson lemma that might be useful
for this task. 

Recently one of the authors with Pott and Reguera proved an apparently
 weak version, featuring a matrix weight, but a scalar Carleson sequence
instead of a matrix sequence and an estimate involving inner products instead
of norms \cite{RPP}. We show that this formulation is optimal in that any improvement of
the statement fails. Indeed, the failure stems from the maximal excentricity
of the matrix weight and not from any increase in the $A_2$ characteristic. We
show that even assuming the $A_2$ condition, one can get an infinite bilinear
embedding. We show that a bilinear estimate can be obtained in terms of the
square root of the conditioning number and this condition is necessary.

\

The failure of BET is natural and requires only a very simple example. It is
however an important notable difference to the scalar case. It is also in a contrast to a positive result on a matrix two--weighted $T1$ theorem 
by Bickel--Culiuc--Treil--Wick \cite{BCTW}. 

There is an array
of difficulties encountered in the task to find various sharp weighted
estimates in the matrix weighted case. Further, it is also very difficult to
get definite negative answers. As of today, most optimal estimates elude us. We mention \cite{IKP2017} for some interesting positive results in the matrix weighted setting.
The first quantitative 
estimate for the Hilbert transform in the matrix weighted setting was given in \cite{BPW2016} as well as 
 upper and lower square function estimates. Some of these estimates in \cite{BPW2016} were close to their scalar weighted analogs, 
 but none of them matched the sharp scalar estimate.
We now know that the estimate for the square function with matrix weight does not
change its dependence on the $A_2$ characteristic as compared to the scalar
case. The result is due to Hyt\"onen--Petermichl--Volberg \cite{HPV2019} and Treil \cite{T} where both use a
stopping time argument known as sparse domination. The classical sharp scalar
weighted result by Hukovi\'c--Treil--Volberg \cite{HTV} was proved via Bellman functions
and only required the simple Carleson Lemma, not the bilinear version. The
estimate for the Hilbert transform is still open, with best to date estimate
by Nazarov--Petermichl--Treil--Volberg \cite{NPTV2017}, missing the sharp
conjecture by a half power of the $A_2$ characteristic. For most known
operator norms, the question of sharpness is unsettled, but there usually is
just a raised power on the dependence of the $A_2$ charateristic, not complete
failure of the estimate, such as what we see in the case of BET.

\

The first proof of a version of the scalar bilinear Carleson lemma is found
in \cite{PW2002} and \cite{P2007} and was rather complicated, using tools and a
construction implicit in the seminal article by Nazarov--Treil--Volberg
\cite{NTV1999}. The argument features a rather cleverly gaged Bellman
function and three conditions on the measure sequence rather than one. It was
understood for some time by the experts that two of the arising conditions
were redundant. We show that this redundancy is still true in the presence of
a matrix weight and a matrix sequence. A previous result only allowed for
scalar Carleson sequences. Since BET fails with matrix Carleson sequences, this redundancy result does not have this
particular application, but it is an interesting estimate in its own right,
useful for other matrix weighted tasks such as certain maximal function estimates. 
It can also be used in combination with CET to alter the testing condition, i.e. changing the Carleson sequence. 
This step was important in some scalar proofs in the non--homogenous setting, see \cite{TTV} and \cite{DP2019}.

\section{Notation and detailed history}

Let us say $Q_0 = [0, 1]$ is endowed with a dyadic filtration and let
$\mathcal{D}$ be the dyadic grid. We call a $d \times d$ matrix--valued
function $W$ a weight if $W (x)$ is positive semidefinite
almost everywhere and if $W$ and $W^{- 1}$ are locally integrable. One defines
$L_{\mathbbm{C}^d}^2 (W)$ to be the set of vector functions with
\begin{eqnarray*}
  \| f \|^2_{L_{\mathbbm{C}^d}^2 (W)} = \int_{Q_0} \| W^{1 / 2} (x) f (x)
  \|_{\mathbbm{C}^d}^2 \mathd x = \int_{Q_0} \langle W (x) f (x), f (x)
  \rangle_{\mathbbm{C}^d} \mathd x < \infty .
\end{eqnarray*}
When $d = 1$ this becomes the classical weighted space $L_{\mathbbm{C}}^2
(w)$. Let us denote by $\langle \cdot \rangle_Q$ the average of a scalar,
vector or matrix function over the cube $Q$. By $\| \cdot \|_{\tmop{op}}$ we
mean the operator norm of the matrix.

\

The dyadic formulation of the Carleson embedding theorem reads as follows:
let $(\alpha_Q)$ be a sequence of non--negative scalars and let $w$ be a
weight. Then for $f$ supported on $Q_0$
\begin{eqnarray*}\lefteqn{
  \frac{1}{| K |} \sum_{Q \in \mathcal{D} (K)} \alpha_Q \langle w \rangle^2_Q
  \lesssim \langle w \rangle_K, \forall K \in \mathcal{D} (Q_0)}\\
 & \Leftrightarrow & \sum_{Q \in \mathcal{D} (Q_0)} \alpha_Q \langle f w^{1/2}
  \rangle^2_Q \lesssim \| f \|^2_{L^2_{\mathbbm{C}}} .
\end{eqnarray*}

There is no difference if the weight is identical to $1$, with the proof in \cite{NT1997} 
being exactly the same when switching to weighted averages and renormalizing appropriately. One can rewrite the
assumption and conclusion renaming $\beta_Q = \alpha_Q \langle w \rangle^2_Q$ and $g=fw^{-1/2}$
with
\begin{eqnarray*}
  \frac{1}{w (K)} \sum_{Q \in \mathcal{D} (K)} \beta_Q \lesssim 1, \forall K
  \in \mathcal{D} (Q_0) \Leftrightarrow \sum_{Q \in \mathcal{D} (Q_0)} \beta_Q
  \langle g \rangle^2_{Q, w} \lesssim \| g \|^2_{L^2_{\mathbbm{C}} (w)},
\end{eqnarray*}
with
$
  \langle g \rangle_{Q, w} = \frac{1}{w (Q)} \int_Q g w.
$
There is as an immediate consequence a CET for matrix Carleson sequences
$(A_Q)$ and the matrix weight $W =\mathbf{1}$:
\begin{eqnarray*}
  \frac{1}{| K |} \sum_{Q \in \mathcal{D} (K)} A_Q \lesssim \mathbf{1},
  \forall K \in \mathcal{D} (Q_0) \Leftrightarrow \sum_{Q \in \mathcal{D}
  (Q_0)} \langle A_Q \langle f \rangle_Q, \langle f \rangle_Q
  \rangle_{\mathbbm{C}^d} \lesssim \| f \|^2_{L^2_{\mathbbm{C}^d}} .
\end{eqnarray*}
Here, the left hand side inequality is understood in the sense of operators
and its necessity `$\Leftarrow$' is seen by testing the right hand inequality
on functions $f = e\mathbf{1}_K$ for any vector $e$ and any $K \in
\mathcal{D} (Q_0)$. The deduction of the implication `$\Rightarrow$' from the
scalar case is easily observed by taking trace and using equivalence of norms
at the cost of a dimensional constant. See \cite{NPTV2002} for details. A weighted
version does not follow from the scalar case, though. The matrix weighted CET
was proved only in 2015:

\begin{theorem}[Culiuc, Treil]
  \label{WCET}Let $W$ be a matrix weight of size $d \times d$ and $(A_Q)$ a
  sequence of positive semidefinite matrices of size $d \times d$. Then for
  $f$ supported in $Q_0$
  \begin{eqnarray*}
  \lefteqn{
    \frac{1}{| K |} \sum_{Q \in \mathcal{D} (K)} \langle W \rangle_Q A_Q
    \langle W \rangle_Q \lesssim \langle W \rangle_K, \forall K \in \mathcal{D}
    (Q_0)}\\
    &\Leftrightarrow& \sum_{Q \in \mathcal{D} (Q_0)} \| A_Q^{1 / 2} \langle W^{1
    / 2} f \rangle_Q \|_{\mathbbm{C}^d}^2 \lesssim \| f
    \|^2_{L^2_{\mathbbm{C}^d}}
      \end{eqnarray*}
\end{theorem}

The successful argument added a twist to the `Bellman function with a
parameter' invented in \cite{NPTV2002}, that managed to make certain
non--commutative obstacles disappear, which arise from differentiating
functions with matrix variables.

Theorem \ref{WCET} implies via a simple linearization argument the norm
estimate below on a Doob type maximal function with matrix measure. 
The obtained norm estimate does not assume the $A_2$ condition. See for example 
\cite{RPP} for an exposition of the well known argument as well as a motivation
of the definition of this maximal function.

\begin{theorem}
  \label{MAX}Let $W$ be a matrix weight of size $d \times d$ then $M_W : L^2_{\mathbbm{C}^d} \rightarrow L^2_{\mathbbm{C}}$ defined by
  \begin{eqnarray*}
    M_W f (x) = \sup_{Q \in \mathcal{D} (Q_0), x \in Q} \| W^{1 / 2} (x)
    \langle W \rangle^{- 1}_Q \langle W^{1 / 2} f \rangle_Q
    \|_{\mathbbm{C}^d},
  \end{eqnarray*}
  is bounded. 
\end{theorem}

We now give some background on the occurrance and motivation for the weighted
bilinear Carleson lemma in the scalar case. In 1999, Nazarov--Treil--Volberg published their
paper \cite{NTV1999} on the necessary and sufficient conditions on the scalar weights $w$ and
$v$ so that the operators $T_{\sigma} : L^2 (w) \rightarrow L^2 (v), h_{_{} Q}
\mapsto \sigma_Q h_Q$ are uniformly bounded. Here
\begin{eqnarray*}
  h_Q = - | Q |^{- 1 / 2} \mathbf{1}_{Q_-} + | Q |^{- 1 / 2}
  \mathbf{1}_{Q_+}
\end{eqnarray*}
are the Haar functions and $\sigma_Q = \pm 1$. Wittwer then proved that
$T_{\sigma} : L^2 (w) \rightarrow L^2 (w)$ has operator norm uniformly bounded
by a linear function of the scalar $A_2$ characteristic of the weight
\begin{eqnarray}
  [w]_{A_2} = \sup_Q \langle w \rangle_Q \langle w^{- 1} \rangle_Q .
  \label{sA2}
\end{eqnarray}
This estimate is sharp, as it was shown to be true for classical operators
such as the Hilbert transform \cite{P2007}. In these early proofs the idea consisted of
`disbalancing' the Haar basis by finding coefficients $a^w_Q$ and $b^w_Q$
\begin{eqnarray}
  h_Q = a^w_Q \mathbf{1}_Q + b^w_Q h^w_Q \label{whaar}
\end{eqnarray}
so that the system $h^w_Q$ becomes orthonormal in $L^2 (w)$ and similarly for
$w^{- 1}$. In its dualized form, the operators $T_{\sigma}$ are uniformly
dominated by the sum
\begin{eqnarray*}
  \sum_{Q \in \mathcal{D} (Q_0)} | \langle w^{1/2}f, h_Q \rangle_{L^2} \langle w^{-1/2}g, h_Q
  \rangle_{L^2} | \lesssim [w]_{A_2} \| f \|_{L_{\mathbbm{C}}^2 } \| g
  \|_{L_{\mathbbm{C}}^2 } .
\end{eqnarray*}
When replacing the Haar functions by the equation (\ref{whaar}) above, we obtain
four sums. The sum featuring the two non--cancellative terms arising from the
indicator functions being the most difficult. It is for this sum that the
bilinear Carleson lemma came to life. The first version appeared in \cite{NTV1999}
in the two weight case with a rather complicated necessary condition on the
two weights $w$ and $v$. A simplified version with easy to test conditions on
the single weight $w$ appeared in \cite{PW2002} \cite{P2007}:
\begin{eqnarray*}\lefteqn{
  \frac{1}{| K |} \sum_{Q \in K} \alpha_Q \lesssim 1, \forall K \in \mathcal{D} (Q_0),}\\
  &&\quad \frac{1}{| K |} \sum_{Q
  \in K} \frac{\alpha_Q}{\langle w \rangle_Q} \lesssim \langle w^{- 1}
  \rangle_K, \forall K \in \mathcal{D} (Q_0), \\
  &&\quad\frac{1}{| K |} \sum_{Q \in K} \frac{\alpha_Q}{\langle w^{- 1}
  \rangle_Q} \lesssim \langle w \rangle_K, \forall K \in \mathcal{D} (Q_0)\\
  &\Rightarrow& \sum_{Q \in \mathcal{D} (Q_0)} \frac{\alpha_Q}{\langle w
  \rangle_Q \langle w^{- 1} \rangle_Q} \langle w^{1 / 2} | f | \rangle_Q
  \langle w^{- 1 / 2} | g | \rangle_Q \lesssim \| f \|_{L^2_{\mathbbm{C}^d}} 
  \| g \|_{L^2_{\mathbbm{C}^d}} \text{} .
\end{eqnarray*}
The conditions are not necessary. Indeed it has been known that two conditions
are redundant:
\begin{eqnarray}\lefteqn{
  \frac{1}{| K |} \sum_{Q \in K} \alpha_Q \lesssim 1, \forall K \in \mathcal{D}(Q_0)\Rightarrow}\\ 
  \nonumber&& \frac{1}{| K
  |} \sum_{Q \in K} \frac{\alpha_Q}{\langle w \rangle_Q} \lesssim \langle w^{-
  1} \rangle_K \tmop{and} \frac{1}{| K |} \sum_{Q \in K}
  \frac{\alpha_Q}{\langle w^{- 1} \rangle_Q} \lesssim \langle w \rangle_K, \forall K \in \mathcal{D}(Q_0).
  \label{red}
\end{eqnarray}

But not even the strongest condition is necessary for the conclusion of the
bilinear embedding. The interest lied in the simpler testing condition and its
successful application to sharp weighted estimates. The first proofs of this
bilinear lemma were quite complicated.

Indeed, in 2008 Nazarov--Treil--Volberg \cite{NTV2008} characterized two--weight estimates for
individual multipliers $T_{\sigma}$ again in the flavor of a $T 1$ theorem.
There is a positive extension to the matrix case by
Bickel--Culiuc--Treil--Wick \cite{BCTW}.

\

Here is the best to date matrix weighted version of BET featuring a scalar
sequence and an inner product:

\begin{theorem}[Petermichl, Pott, Reguera]
  \label{siBET}Let $(\alpha_Q)$ be a sequence of non--negative scalars. Then
  for $f, g$ supported in $Q_0$
  \begin{eqnarray*}\lefteqn{
    \frac{1}{| K |} \sum_{Q \in \mathcal{D} (K)} \alpha_Q \lesssim 1, \forall K
    \in \mathcal{D} (Q_0)}\\
    &\Rightarrow & 
    \sum_{Q \in \mathcal{D} (Q_0)} \alpha_Q | \langle \langle W
    \rangle^{- 1}_Q \langle W^{1 / 2} f \rangle_Q, \langle W^{- 1} \rangle^{-
    1}_Q \langle W^{- 1 / 2} g \rangle_Q \rangle_{\mathbbm{C}^d} | \lesssim \|
    f \|_{L^2_{\mathbbm{C}^d}}  \| g \|_{L^2_{\mathbbm{C}^d}}.
  \end{eqnarray*}
\end{theorem}

We also study matrix analogs for the redundancy whose scalar version is the
implication in (\ref{red}). Here is the redundancy with scalar sequence and
matrix weight, which was proved recently:

\begin{theorem}[Petermichl, Pott, Reguera]
  \label{sRED}Let $(\alpha_Q)$ be a non--negative sequence and $W$ a matrix
  weight, then
  \begin{eqnarray*}\lefteqn{
    \frac{1}{| K |} \sum_{Q \in \mathcal{D} (K)} \alpha_Q \leqslant 1, \forall
    K \in \mathcal{D} (Q_0)}\\
   & \Rightarrow & \frac{1}{| K |} \sum_{Q \in \mathcal{D} (K)} \alpha_Q \langle
    W^{- 1} \rangle_Q^{- 1} \leqslant 4 \langle W \rangle_K, \forall K \in
    \mathcal{D} (Q_0) .
  \end{eqnarray*}
\end{theorem}

\section{Main Results}

In this paper we prove the sharpness of the bilinear Embedding Theorem
\ref{siBET}, i.e. the failure of any improvement:

\begin{theorem}
  \label{BETfailure}Let $(A_Q)$ be a sequence of positive semidefinite
  matrices. The Carleson condition
  \begin{eqnarray}\label{conditionBETfailure}
    \frac{1}{| K |} \sum_{Q \in \mathcal{D} (K)} A_Q \lesssim \mathbf{1},
    \forall K \in \mathcal{D} (Q_0)
  \end{eqnarray}
  does not imply the existence of a constant so that
  \begin{eqnarray}\label{miBETfailure}\lefteqn{
    \sum_{Q \in \mathcal{D} (Q_0)} |  \langle A_Q \langle W \rangle^{- 1}_Q
    \langle W^{1/2} f \rangle_Q, \langle W^{- 1} \rangle^{- 1}_Q \langle W^{- 1/2} g
    \rangle_Q \rangle_{\mathbbm{C}^d} |}\\
    \nonumber& \lesssim  &\| f \|_{L_{\mathbbm{C}^d}^2
    } \| g \|_{L_{\mathbbm{C}^d}^2}
  \end{eqnarray}
  for all matrix weights $W$ and functions $f, g$ supported in $Q_0$. The
  Carleson condition also does not imply
  \begin{eqnarray}\label{mnBETfailure}\lefteqn{
    \sum_{Q \in \mathcal{D} (Q_0)}  \| A^{1 / 2}_Q \langle W \rangle^{- 1}_Q
    \langle W^{1/2} f \rangle_Q \|_{\mathbbm{C}^d} \| A^{1 / 2}_Q \langle W^{- 1}
    \rangle^{- 1}_Q \langle W^{- 1/2} g \rangle_Q \|_{\mathbbm{C}^d}}\\
    \nonumber & \lesssim & \|
    f \|_{L_{\mathbbm{C}^d}^2} \| g \|_{L_{\mathbbm{C}^d}^2 }
  \end{eqnarray}
  for all matrix weights $W$ and functions $f, g$ supported in $Q_0$, not even
  if we only allow scalar sequences $A_Q = m_Q$. 
\end{theorem}

We find that a natural bilinear Carleson embedding theorem only holds under a
very strong condition on the matrix weight:

\begin{definition}
  Let $W$ be a matrix weight. We define the conditioning number of a matrix
  weight as
  \begin{eqnarray*}
    [W]_{C_2} = \sup_x \kappa (W (x)) = \sup_x \| W (x) \|_{\tmop{op}} \| W^{-
    1} (x) \|_{\tmop{op}} = \frac{\lambda_{\max} (W (x))}{\lambda_{\min} (W
    (x))},
  \end{eqnarray*}
  where $\lambda_{\max}$ and $\lambda_{\min}$ return the maximal respectively
  minimal eigenvalue.
\end{definition}

This is in a sharp contrast to the matrix $A_2$ condition in \cite{TV1997}, the analog of the scalar version in equation (\ref{sA2}):

\begin{definition}
  Let $W$ be a matrix weight. Then the matrix $A_2$ condition is
  \begin{eqnarray*}
    [W]_{A_2} = \sup_Q \| \langle W \rangle^{1 / 2}_Q \langle W^{- 1}
    \rangle^{1 / 2}_Q \|_{\tmop{op}}^2 .
  \end{eqnarray*}
\end{definition}

See \cite{TV1997}  for a list of elementary properties of this
characteristic, such as $[W]_{A_2} \geqslant 1$. We prove:

\begin{theorem}
  \label{C2BET}Let $(A_Q)$ be a sequence of positive semidefinite matrices
  and $W$ a matrix weight. Then
  \begin{eqnarray*}\lefteqn{
    \frac{1}{| K |} \sum_{Q \in \mathcal{D} (K)} A_Q \lesssim \mathbf{1},
    \forall K \in \mathcal{D} (Q_0)}\\
      &  \Rightarrow & \sum_{Q \in \mathcal{D} (Q_0)}  \| A^{1 / 2}_Q \langle W
    \rangle^{- 1}_Q \langle W^{1 / 2} f \rangle_Q \|_{\mathbbm{C}^d} \| A^{1 /
    2}_Q \langle W^{- 1} \rangle^{- 1}_Q \langle W^{- 1 / 2} g \rangle_Q
    \|_{\mathbbm{C}^d} \\
    && \quad \lesssim [W]^{1 / 2}_{C_2} \| f
    \|_{L_{\mathbbm{C}^d}^2} \| g \|_{L_{\mathbbm{C}^d}^2 .}
  \end{eqnarray*}

  The condition $[W]_{C_2}$ is necessary and the power $1 / 2$ is optimal.
\end{theorem}

Further, we deduce from Theorem \ref{sRED} in full generality:

\begin{theorem}
  \label{RED}Let $(B_Q)$ be a sequence of positive semidefinite matrices such
  that
  \begin{eqnarray*}
    \frac{1}{| {L} |} \sum_{Q \in \mathcal{D} ({
    L})} B_Q \lesssim \mathbf{1}, \forall { L} \in \mathcal{D}
    (Q_0) .
  \end{eqnarray*}
  Then 
  \begin{eqnarray*}\lefteqn{
    \frac{1}{| { K} |} \sum_{Q \in \mathcal{D} ({
    K})} \langle B_Q \langle W \rangle_{{ K}}^{- 1 / 2} \langle
    W^{- 1} \rangle^{- 1 / 2}_Q { e}, \langle W \rangle^{- 1 /
    2}_{{ K}} \langle W^{- 1} \rangle^{- 1 / 2}_Q {
    e} \rangle_{\mathbbm{C}^d}}\\
    &&\quad \lesssim  \| { e}
    \|_{\mathbbm{C}^d}^2, \forall {K} \in \mathcal{D} (Q_0) .
  \end{eqnarray*}
  and
  \begin{eqnarray*}\lefteqn{
    \frac{1}{| { K} |} \sum_{Q \in \mathcal{D} ({
    K})} \langle B_Q \langle W^{- 1} \rangle^{- 1 / 2}_Q \langle W \rangle^{-
    1 / 2}_{{ K}} { e}, \langle W^{- 1} \rangle^{- 1
    / 2}_Q \langle W \rangle^{- 1 / 2}_{{K}} {e}
    \rangle_{\mathbbm{C}^d} }\\
    &&\quad \lesssim \| { e} \|_{\mathbbm{C}^d}^2,
    \forall {\color{black} K} \in \mathcal{D} (Q_0) .
  \end{eqnarray*}
  for all vectors ${e
  }$ and all matrix weights $W$.
  \end{theorem}

Notice that in Theorem \ref{RED} putting ${ f} = \langle W \rangle^{- 1 / 2}_{{ K}}
{ e}$ the last inequality becomes
\begin{eqnarray*}\lefteqn{
  \frac{1}{| {\color{black} K} |} \sum_{Q \in \mathcal{D} ({K})}
  \langle B_Q \langle W^{- 1} \rangle^{- 1 / 2}_Q { f}, \langle
  W^{- 1} \rangle^{- 1 / 2}_Q {\color{black} f} \rangle_{\mathbbm{C}^d}}\\
  &&\quad
  \lesssim \| \langle W \rangle^{1 / 2}_{{ K}} {
  { f}} \|_{\mathbbm{C}^d}^2, \forall { K} \in
  \mathcal{D} (Q_0),
\end{eqnarray*}
that rewrites as the operator inequality resembling (\ref{red}):
\begin{eqnarray*}
  \frac{1}{| { K} |} \sum_{Q \in \mathcal{D} ( K)}
  \langle W^{- 1} \rangle^{- 1 / 2}_Q B_Q \langle W^{- 1} \rangle^{- 1 / 2}_Q
  \lesssim \langle W \rangle_K,\forall {K} \in \mathcal{D} (Q_0) .
\end{eqnarray*}

\section{Bilinear embedding theorem}

We prove Theorems \ref{BETfailure} and \ref{C2BET} in this section.

\begin{proof}[Theorem \ref{C2BET}]
  First, we can assume that the sequence consists of scalars, by switching to
  the Carleson sequence of maximal eigenvalues $\alpha_Q = \| A_Q
  \|_{\tmop{op}}$ at the loss of a dimensional constant. Then we have to show
  the inequality
  \begin{eqnarray}\lefteqn{
    \sum_{Q \in \mathcal{D} (Q_0)} \alpha_Q  \| \langle W \rangle^{- 1}_Q
    \langle W^{1 / 2} f \rangle_Q \|_{\mathbbm{C}^d} \| \langle W^{- 1}
    \rangle^{- 1}_Q \langle W^{- 1 / 2} g \rangle_Q \|_{\mathbbm{C}^d}}\\
  \nonumber&&\quad  \lesssim [W]^{1 / 2}_{C_2} \| f \|_{L_{\mathbbm{C}^d}^2} \| g
    \|_{L_{\mathbbm{C}^d}^2} . \label{C2ineq}
  \end{eqnarray}

  In order to prove this inequality, let us define
  \begin{eqnarray*}
    {\mu} (\mathcal{K}) = \sum_{Q \in \mathcal{K}} \alpha_Q
  \end{eqnarray*}
  for any collection $\mathcal{K}$ of dyadic cubes. Let $F$ be any
  non--negative function defined on the dyadic cubes. Then let denote $\{ F
  (Q) > \lambda \}$ the collection of cubes so that $F (Q) > \lambda$. It
  follows that
  \begin{eqnarray*}
    \int^{\infty}_0 \left( \sum_{Q \in \{ F (Q) > \lambda \}} \alpha_Q \right)
    \mathd \lambda = \int^{\infty}_0 {\mu} (\{ F (Q) > \lambda \}) \mathd
    \lambda = \sum_{Q \in \mathcal{D} (Q_0)} F (Q) \alpha_Q,
  \end{eqnarray*}
  which is the classical fact on Choquet integrals. Let us pose
  \begin{eqnarray*}
    F (Q) = \| \langle W \rangle^{- 1}_Q \langle W^{1 / 2} f \rangle_Q
    \|_{\mathbbm{C}^d} \| \langle W^{- 1} \rangle^{- 1}_Q \langle W^{- 1 / 2}
    g \rangle_Q \|_{\mathbbm{C}^d}
  \end{eqnarray*}
  and let for any $\lambda > 0$ the set $\mathcal{J}_{\lambda}$ denote the
  collection of maximal dyadic intervals for which $F (Q) > \lambda$. So the
  integrand above becomes
  \begin{eqnarray*}
    \sum_{Q \in \{ F (Q) > \lambda \}} \alpha_Q \leqslant \sum_{Q \in
    \mathcal{J}_{\lambda}} \sum_{Q' \in \mathcal{D} (Q)} \alpha_Q \leqslant
    \sum_{Q \in \mathcal{J}_{\lambda}} | Q | .
  \end{eqnarray*}
  Now let
  \begin{eqnarray*}
    \Phi (x) = M_W f (x) M_{W^{- 1}} g (x) .
  \end{eqnarray*}
  Observe that with
  \begin{eqnarray*}
    A^W_Q f (x) \assign W^{1 / 2} (x) \langle W \rangle^{- 1}_Q \langle W^{1 /
    2} f \rangle_Q 
  \end{eqnarray*}
  and
  \begin{eqnarray*}
    A^{W^{- 1}}_Q g (x) = W^{- 1 / 2} (x) \langle W^{- 1} \rangle^{- 1}_Q
    \langle W^{- 1 / 2} g \rangle_Q
  \end{eqnarray*}
  by Cauchy Schwarz for all $x$ and all $Q$ there holds
  \begin{eqnarray*}\lefteqn{
    F (Q) }\\
    & = & \| W^{- 1 / 2} (x) W^{1 / 2} (x) \langle W \rangle^{- 1}_Q
    \langle W^{1 / 2} f \rangle_Q \|_{\mathbbm{C}^d} \\
    &&\quad \| W^{1 / 2} (x) W^{- 1 /
    2} (x) \langle W^{- 1} \rangle^{- 1}_Q \langle W^{- 1 / 2} g \rangle_Q
    \|_{\mathbbm{C}^d}\\
    & \leqslant & \| W^{- 1 / 2} (x) \|_{\tmop{op}} \| W^{1 / 2} (x)
    \|_{\tmop{op}} \| A^W_Q f (x) \|_{\mathbbm{C}^d} \| A^{W^{- 1}}_Q g (x)
    \|_{\mathbbm{C}^d} .
  \end{eqnarray*}
  Observe that
  \begin{eqnarray*}
    \sup_y \| W^{- 1 / 2} (y) \|_{\tmop{op}} \| W^{1 / 2} (y) \|_{\tmop{op}} =
    [W]^{1 / 2}_{C_2}
  \end{eqnarray*}
  and for $x \in Q$ we get
  \begin{eqnarray*}
    F (Q) \lesssim [W]^{1 / 2}_{C_2} \Phi (x) .
  \end{eqnarray*}
  So if $x \in Q$ with $F (Q) > \lambda$ then also $[W]^{1 / 2}_{C_2} \Phi (x)
  > \lambda$. So
  \begin{eqnarray*}
    \sum_{Q \in \mathcal{J}_{\lambda}} | Q | \leqslant | \{ x \in Q_0 : [W]^{1
    / 2}_{C_2} \Phi (x) > \lambda \} | .
  \end{eqnarray*}
  Integrating with respect to $\mathd \lambda$ gives
  \begin{eqnarray*}\lefteqn{
    \sum_{Q \in \mathcal{D} (Q_0)} \alpha_Q F (Q)}\\
    & \leqslant & \int^{\infty}_0 |
    \{ x \in Q_0 : [W]^{1 / 2}_{C_2} \Phi (x) > \lambda \} | \mathd \lambda \\
    &=&
    [W]^{1 / 2}_{C_2} \int_{Q_0} M_W f (x) M_{W^{- 1}} g (x) \mathd x.
  \end{eqnarray*}
  Using the above estimate of the maximal function, Theorem \ref{MAX} and an
  application of Cauchy Schwarz finishes the proof of estimate (\ref{C2ineq}).
\end{proof}

Now, we prove Theorem \ref{BETfailure}, the failure of any improvement of
Theorem \ref{siBET}.

\begin{proof}[Theorem \ref{BETfailure}]
  Let us take the case $d = 2$ and let $a$ and $b$ be orthogonal unit vectors.
  Let $W = a a^{\ast} + \varepsilon^2 b b^{\ast}$ and thus $W^{- 1} = a
  a^{\ast} + \varepsilon^{- 2} b b^{\ast}$. Letting $f = W^{1 / 2}
  b\mathbf{1}_{Q_0}$ and $g = W^{- 1 / 2} a\mathbf{1}_{Q_0}$ we get 
  \begin{eqnarray*}
    \| f \|_{L_{\mathbbm{C}^d}^2} = \langle W b, b \rangle_{\mathbbm{C}^d}^{1 / 2} =
    \varepsilon \; \tmop{ and } \; \| g \|_{L_{\mathbbm{C}^d}^2} = \langle W^{- 1} a,
    a \rangle_{\mathbbm{C}^d}^{1 / 2} = 1.
  \end{eqnarray*}
  This gives us the order $\varepsilon$ for the right hand sides of inequalities (\ref{miBETfailure}) and (\ref{mnBETfailure}).
  Further, we choose the Carleson sequence $A_{Q_0} =\mathbf{1}$ and $A_Q=\mathbf{0}$ for all other cubes. 
  The Carleson intensity in inequality (\ref{conditionBETfailure}) is 1. We get for the sum in inequality (\ref{mnBETfailure}) only one term:
  \begin{eqnarray}
    \| A^{1 / 2}_{Q_0} \langle W \rangle^{- 1}_{Q_0} \langle W^{1 / 2} f
    \rangle_{Q_0} \|_{\mathbbm{C}^d} \| A^{1 / 2}_{Q_0} \langle W^{- 1}
    \rangle^{- 1}_{Q_0} \langle W^{- 1 / 2} g \rangle_{Q_0} \|_{\mathbbm{C}^d}
    = 1. \label{BETtest}
  \end{eqnarray}
  Letting $\varepsilon \to 0$, we have shown the failure of conclusion (\ref{mnBETfailure}).
  The conditioning number of the weight $W$ is $[W]_{C_2} = \varepsilon^2$.
  Again letting $\varepsilon \rightarrow 0$
  shows the necessity of $[W]^{1 / 2}_{C_2}$ occurring in Theorem \ref{C2BET}. 
  To see
  that the inequality (\ref{miBETfailure}) also fails, choose $A_{Q_0} =
  2^{- 1} (a+b) (a+b)^{\ast} $ and $A_Q= \mathbf{0}$ for all other cubes. We see that 
  \begin{eqnarray*}
    \langle A^{1 / 2}_{Q_0} \langle W \rangle^{- 1}_{Q_0} \langle W^{1 / 2} f
    \rangle_{Q_0}, A^{1 / 2}_{Q_0} \langle W^{- 1} \rangle^{- 1}_{Q_0} \langle
    W^{- 1 / 2} g \rangle_{Q_0} \rangle_{\mathbbm{C}^d} \gtrsim 1,
  \end{eqnarray*}
  showing the failure of conclusion (\ref{miBETfailure}).
  It is easy to replace the matrix Carleson sequence by the scalar sequence
  $\alpha_{Q_0} = 1$ and $\alpha_Q=0$ for all other cubes as the expression
  (\ref{conditionBETfailure}) does not change.
\end{proof}

As we see, the bilinear Carleson Lemma fails violently - using constant
weights with a high discrepancy in their eigenvalues is sufficient. It is
natural that this quantity determines the growth of the estimate and not
features measured by the matrix $A_2$ characteristic. Notice that $[W]_{A_2} =
1$ for our example. \

\section{Redundant Carleson condition}

In this section we prove Theorem \ref{RED}. We sketch the proof of Theorem
\ref{sRED}. Consider the matrix valued Bellman function of matrix variables
$U, V$ and scalar variable $m$
\begin{eqnarray*}
  B (U, V, m) = U - (m + 1)^{- 1} V^{- 1} .
\end{eqnarray*}
This function has domain $\mathbf{1} \leqslant V^{1 / 2} U V^{1 / 2}$ and $0
\leqslant m \leqslant 1.$ There holds the size estimate
\begin{eqnarray*}
  0 \leqslant B (U, V, m) \leqslant U.
\end{eqnarray*}
Indeed, $0 \leqslant U - V^{- 1} \leqslant U - (m + 1)^{- 1} V^{- 1} \leqslant
U.$ The function $B$ is also concave: Dropping the linear dependence on $U$,
its Hessian acting on the matrix difference $\Delta V$ and scalar $\Delta m$
is a positive multiple of
\begin{eqnarray*}\lefteqn{
  - 2 V^{- 1} \Delta V V^{- 1} \Delta V V^{- 1} }\\
  &&\quad - 2 V^{- 1}
  \Delta V V^{- 1} (m + 1)^{- 1} \Delta m - 2 (m + 1)^{- 2} V^{- 1} (\Delta
  m)^2 .
\end{eqnarray*}
Observe that
\begin{eqnarray*}
  V^{- 1} \Delta V V^{- 1} \Delta V V^{- 1}  \geqslant 0, \quad
  V^{- 1} (\Delta m)^2 \geqslant 0.
\end{eqnarray*}
If we add positive multiples of these non--negative terms to the Hessian, we can write it as a
negative perfect square and therefore concavity follows. We also have
\begin{eqnarray*}
  \partial B / \partial m = (m + 1)^{- 2} V^{- 1} \geqslant 4^{- 1} V^{- 1} .
\end{eqnarray*}

One can check that the variables $$U_K = \langle W \rangle_K, \; V_K = \langle W^{-
1} \rangle_K, \; m_K = | K |^{- 1} \sum_{Q \in \mathcal{D} (K)} \alpha_Q$$ lie
in the domain of $B$. We see that $m_K - | K |^{- 1} \alpha_K = 2^{- 1}
(m_{K_-} + m_{K_+})$. The usual Bellman dynamics argument gives the estimate
on the operator sum:
\begin{eqnarray*}\lefteqn{
  | K | \langle W \rangle_K = | K | U_K
   \geqslant  | K | B (U_K, V_K, m_K)}\\
  && \quad =  | K | (B (U_K, V_K, m_K) \\
  && \quad \quad - B (U_K, V_K, m_K - | K |^{- 1} \alpha_K) +
  B (U_K, V_K, m_K - | K |^{- 1} \alpha_K)\\
  && \quad \geqslant  4^{- 1} V_K^{- 1} \alpha_K + | K | B (U_K, V_K, m_K - | K
  |^{- 1} \alpha_K)\\
  && \quad \geqslant  4^{- 1} V_K^{- 1} \alpha_K + | K_- | B (U_{K_-}, V_{K_-},
  m_{K_-}) + | K_+ | B (U_{K_+}, V_{K_+}, m_{K_+})
\end{eqnarray*}
Iterating this argument gives the desired estimate for scalar sequences
$(\alpha_Q)$. We note that if the Carleson measure is not scalar, then the
function to consider may be
\begin{eqnarray*}
  U - V^{- 1 / 2} (M +\mathbf{1})^{- 1} V^{- 1 / 2} .
\end{eqnarray*}
The concavity of this function is unclear to us. Now we argue that we can
conclude anyways.

\begin{proof}[Theorem \ref{RED}] The required observation is the following:
\begin{eqnarray*}\lefteqn{
  \frac{1}{| L|} \sum_{Q \in \mathcal{D}
  (L)} B_Q \lesssim \mathbf{1}, \forall L\in
  \mathcal{D} (Q_0) }\\
  & \Leftrightarrow & \frac{1}{| L |}
  \sum_{Q \in \mathcal{D} ( L)} \tmop{tr} (B_Q) \lesssim 1,
  \forall { L} \in \mathcal{D} (Q_0)\\
  & \Leftrightarrow & \frac{1}{| { L} |} \sum_{Q \in \mathcal{D}
  ({ L})} \| B_Q \|_{\tmop{op}} \lesssim 1, \forall {
  L} \in \mathcal{D} (Q_0) .
\end{eqnarray*}
The implied constants may depend upon $d$. 
Applying Theorem \ref{sRED} to the Carleson sequence $b_Q = \| B_Q
\|_{\tmop{op}}$ gives us
\begin{eqnarray*}
  \frac{1}{| K |} \sum_{Q \in \mathcal{D} (K)} \langle W^{- 1} \rangle^{- 1 /
  2}_Q b_Q \langle W^{- 1} \rangle_Q^{- 1 / 2} \lesssim \langle W \rangle_K,
  \forall K \in \mathcal{D} (Q_0)
\end{eqnarray*}
and thus since $b_Q \mathbf{1} \geqslant B_Q$ we obtain
\begin{eqnarray*}
  \frac{1}{| K |} \sum_{Q \in \mathcal{D} (K)} \langle W^{- 1} \rangle^{- 1 /
  2}_Q B_Q \langle W^{- 1} \rangle_Q^{- 1 / 2} \lesssim \langle W \rangle_K,
  \forall K \in \mathcal{D} (Q_0)
\end{eqnarray*}
from which it follows  that
\begin{eqnarray*}\lefteqn{
  \frac{1}{| K |} \sum_{Q \in \mathcal{D} (K)} \langle W \rangle^{- 1 /
  2}_K  \langle W^{- 1} \rangle^{- 1 / 2}_Q B_Q \langle W^{- 1} \rangle^{-
  1 / 2}_Q \langle W \rangle^{- 1 / 2}_K}\\ 
  &&\quad \lesssim  \mathbf{1}, \forall K \in
  \mathcal{D} (Q_0) .
\end{eqnarray*}
Also
\begin{eqnarray*}
  \frac{1}{| K |} \sum_{Q \in \mathcal{D} (K)} \langle W^{- 1} \rangle^{- 1 /
  2}_Q b_Q \langle W^{- 1} \rangle_Q^{- 1 / 2} \lesssim \langle W \rangle_K,
  \forall K \in \mathcal{D} (Q_0)
\end{eqnarray*}
implies in particular
\begin{eqnarray*}
  \frac{1}{| K |} \sum_{Q \in \mathcal{D} (K)} \tmop{tr} (\langle W \rangle^{-
  1 / 2}_K \langle W^{- 1} \rangle^{- 1 / 2}_Q b_Q \langle W^{- 1}
  \rangle^{- 1 / 2}_Q \langle W \rangle^{- 1 / 2}_K) \lesssim 1, \forall K \in
  \mathcal{D} (Q_0)
\end{eqnarray*}
and with
\begin{eqnarray*}\lefteqn{
  \tmop{tr} (\langle W \rangle^{- 1 / 2}_K  \langle W^{- 1} \rangle^{- 1 /
  2}_Q b_Q \langle W^{- 1} \rangle^{- 1 / 2}_Q \langle W \rangle^{- 1 / 2}_K)}\\
  &&\quad = \tmop{tr} (^{}  \langle W^{- 1} \rangle^{- 1 / 2}_Q \langle W \rangle^{- 1
  / 2}_K b_Q \langle W \rangle^{- 1 / 2}_K \langle W^{- 1} \rangle^{- 1 /
  2}_Q)
\end{eqnarray*}
thus
\begin{eqnarray*}
  \frac{1}{| K |} \sum_{Q \in \mathcal{D} (K)} \tmop{tr} (\langle W^{- 1}
  \rangle^{- 1 / 2}_Q \langle W \rangle^{- 1 / 2}_K b_Q \langle W
  \rangle^{- 1 / 2}_K \langle W^{- 1} \rangle^{- 1 / 2}_Q) \lesssim 1, \forall
  K \in \mathcal{D} (Q_0)
\end{eqnarray*}
and so
\begin{eqnarray*}
  \frac{1}{| K |} \sum_{Q \in \mathcal{D} (K)} \langle W^{- 1} \rangle^{- 1 /
  2}_Q \langle W \rangle^{- 1 / 2}_K b_Q \langle W \rangle^{- 1 / 2}_K
  \langle W^{- 1} \rangle^{- 1 / 2}_Q \lesssim \mathbf{1}, \forall K \in
  \mathcal{D} (Q_0) .
\end{eqnarray*}
So we also obtain with $b_Q \mathbf{1} \geqslant B_Q$ that
\begin{eqnarray*}
  \frac{1}{| K |} \sum_{Q \in \mathcal{D} (K)} \langle W^{- 1} \rangle^{- 1 /
  2}_Q \langle W \rangle^{- 1 / 2}_K B_Q \langle W \rangle^{- 1 / 2}_K
  \langle W^{- 1} \rangle^{- 1 / 2}_Q \lesssim \mathbf{1}, \forall K \in
  \mathcal{D} (Q_0) .
\end{eqnarray*}

\end{proof}

\

\

\

\


\

\

\

\

\end{document}